\documentclass[12pt]{article}
\usepackage{bbm}
\usepackage{mathrsfs}
\usepackage{amsfonts}
\usepackage{amsmath}
\usepackage{amssymb,amsmath}
\def\cl{\centerline}
\def\vs{\vspace*}

\def\QED{\hfill$\Box$}

\numberwithin{equation}{section}
\newtheorem{theo}{Theorem}[section]
\newtheorem{defi}[theo]{Definition}
\newtheorem{exam}[theo]{Example}
\newtheorem{coro}[theo]{Corollary}

\newtheorem{prop}[theo]{Proposition}

\begin{document}

\cl{{\large\bf Hom-Novikov Algebras and Hom-Novikov-Poisson Algebras
} \footnote{Corresponding author: lmyuan@mail.ustc.edu.cn}}
\vspace{2pt}

\cl{ Lamei Yuan$^{\,\ddag}$, Hong You$^{\,\ddag\,\dag}$ }
 \cl{\small{ $\ddag$Academy of Fundamental and Interdisciplinary
 Sciences,}}
\cl{\small{Harbin Institute of Technology, Harbin 150080, China}}
\cl{\small{$^{\dag}$Department of Mathematics, Suzhou University,
Suzhou 200092, China}} \cl{\small E-mail: lmyuan@mail.ustc.edu.cn,
youhong@suda.edu.cn
 }\vs{6pt}

{\small\parskip .005 truein \baselineskip 3pt \lineskip 3pt

\noindent{\bf Abstract.} The purpose of this paper is to study
Hom-Novikov algebras and Hom-Novikov-Poisson algebras, both of which
were defined by Yau. In the paper, we give several constructions
leading us to some interesting examples of Hom-Novikov algebras and
Hom-Novikov-Poisson algebras. Also, we introduce the notion of
quadratic Hom-Novikov algebras and provide some properties.
 \vs{5pt}


\parskip .001 truein\baselineskip 6pt \lineskip 6pt

\vs{10pt}

\cl{\bf\S1. \
Introduction}\setcounter{section}{1}\setcounter{equation}{0}
\vs{6pt}

Novikov algebras were introduced in connection with the Poisson
brackets of hydrodynamic type \cite{BN,DN1,DN2} and Hamiltonian
operators in the formal variational calculus \cite{GD,Xu1}. The
theoretical study of Novikov algebras was started by Zel'manov
\cite{Z} and Filipov \cite{F}. But the term ``Novikov algebra" was
first used by Osborn \cite{O1}.  The left multiplication operators
of a Novikov algebra form a Lie algebra. Thus, it is effective to
relate the study of Novikov algebras to the theory of Lie algebras
\cite{BM1}. Novikov algebras are a special class of left-symmetric
algebras (or under other names such as pre-Lie algebras,
quasi-associative algebras and Vinberg algebras), arising from the
study of affine manifolds, affine structures and convex homogeneous
cones \cite{A,K,V}. Left-symmetric algebras have close relations
with many important fields in mathematics and mathematical physics,
such as infinite-dimensional Lie algebras \cite{BN,DN2}, classical
and quantum Yang-Baxter equation \cite{ES,GS}, quantum field theory
\cite{CK} and so on.

 Novikov-Poisson
algebras were originally introduced by Xu \cite{Xu1, Xu2}  motivated
by studying simple Novikov algebras and irreducible modules. Every
Novikov-Poisson algebra with associative commutative unity can be
constructed from an associative commutative derivation algebra. That
there exists a relationship between associative commutative
derivation algebras and Novikov algebras was also pointed out in
\cite{GD}. Infinite-dimensional Novikov algebras over an
algebraically closed field of characteristic $0$ were classified in
\cite{Xu3} under some restrictions on the structure of right and
left multiplication operators. It turned out that such Novikov
algebras are Novikov-Poisson algebras and are obtained from
associative commutative differentially simple algebras with unity.
While it was  proved in \cite{ZT} that a Novikov algebra is simple
over a field of characteristic not $2$ if and only if its
associative commutative derivation algebra is differentially simple.

A quadratic Novikov algebra, introduced in \cite{Bo, ZC}, is a
Novikov algebra with  a symmetric nondegenerate invariant bilinear
form. The motivation for studying quadratic Novikov algebras came
from the fact that Lie algebras or associative algebras with
symmetric nondegenerate  invariant bilinear forms have important
applications in several areas of mathematics and physics, such as
the structure theory of finite-dimensional semisimple Lie algebras,
the theory of complete integrable Hamiltonian systems and the
classification of statistical models over two-dimensional graphs.

Hom-Lie algebras were initially introduced by Hartwig, Larsson and
Silvestrov in \cite{HLS} motivated by examples of deformed Lie
algebras coming from twisted discretizations of vector fields. This
kind of algebras includes Lie algebras as a subclass in which the
twisting map is the identity map. As generalizations of Hom-Lie
algebras, quasi-hom-Lie algebras \cite{LS}, Hom-Lie superalgebras
\cite{AM} and Hom-Lie color algebras \cite{Yuan} were introduced,
respectively. Both the Hom-Lie superalgebras and Hom-Lie color
algebras are special cases of the quasi-hom-Lie algebras. There has
been a lot of progress in the study of Hom-Lie algebras (see e.g.,
\cite{Yau1,Yau2,Yau3,Yau4}).

Hom-associative algebras generalizing associative algebras to a
situation where associativity law is twisted by a linear map were
introduced in \cite{MS}. It is well known that there is always a Lie
algebra associated to an associative algebra via the commutator
bracket. An analogous result holds for Hom-associative algebras,
i.e., the commutator bracket multiplication defined by the
multiplication in a Hom-associative algebra leads naturally to
Hom-Lie algebras. This result was also extended to the Hom-Lie
superalgebra \cite{AM} and the Hom-Lie color algebra \cite{Yuan}
cases.

Following the patterns of Hom-Lie algebras and Hom-associative
algebras, Yau in \cite{Yau} introduced Hom-Novikov algebras, in
which the two defining identities are twisted by a linear map. It
turned out that Hom-Novikov algebras can be constructed from Novikov
algebras, commutative Hom-associative algebras and Hom-Lie algebras
along with some suitable linear maps. Later, Yau in \cite{Yau5}
defined a Hom-Novikov-Poisson algebra as a twisted generalization of
Novikov-Poisson algebras.

In the present paper, we consider both Hom-Novikov algebras and
Hom-Novikov-Poisson algebras. In Section 2, we summarize some
definitions including Novikov algebras, Novikov-Poisson algebras,
Hom-associative algebras, Hom-Lie algebras, Hom-Novikov algebras and
Hom-Novikov-Poisson algebras and related results.

In Section 3, we focus on Hom-Novikov algebras. We establish the
relationships between Hom-Novikov algebras and Hom-Lie algebras as
well as the relations between Hom-Novikov algebras and Novikov
algebras. Also, we provide a construction of Hom-Novikov algebras
from a commutative Hom-associative algebra together with a suitable
linear map but not a derivation.

In Section 4, we extend the notion of quadratic Novikov algebras to
quadratic Hom-Novikov algebras, which are Hom-Novikov algebras with
a symmetric nondegenerate  bilinear form satisfying equation
(\ref{quadratic}). We reduce the case where the twisting map is an
automorphism or an involution to the study of the relationships
between quadratic Hom-Novikov algebras and quadratic Novikov
algebras as well as the relations between quadratic Hom-Novikov
algebras and quadratic Hom-Lie algebras.

In Section 5, we study Hom-Novikov-Poisson algebras. We prove that
the tensor of two Hom-Novikov-Poisson algebras is still a
Hom-Novikov-Poisson algebra. In addition, we give a way to construct
Hom-Novikov-Poisson algebras by using a commutative Hom-associative
algebra along with a derivation, generalizing a construction in
Novikov algebra case due to Dorfman and Gel'fand, in Novikov-Poisson
algebra case due to Xu and in Hom-Novikov algebra case due to Yau.
Finally, we provide some interesting examples.

Throughout this paper, $\mathbb{F}$ denotes a  field of
characteristic zero. All vector spaces and tensor products are
assumed to be over $\mathbb{F}$, except where otherwise indicated.

\vs{8pt}\

\cl{\bf\S2. \ Preliminaries
}\setcounter{section}{2}\setcounter{equation}{0}\setcounter{theo}{0}
\vs{8pt}\

In this section, we summarize some definitions concerning Novikov
algebras and Hom-algebras, and related results. For detailed
discussions and examples we refer the reader to the literatures, for
instance, \cite{HLS,MS,Xu1, Xu2, Yau,Yau5} and references therein.
\begin{defi}{\rm
A Novikov algebra $(\mathcal {A},\mu)$ is a vector space  with a
bilinear product
$\mu:\mathcal{A}\times\mathcal{A}\longrightarrow\mathcal{A}$
satisfying
\begin{eqnarray}
(xy)z &=&(xz)y,\label{Nov2}\\
(xy)z-x(yz) &=& (yx)z -y(xz),\label{Nov1}
\end{eqnarray}
for all $x,y,z\in \mathcal {A}$, and where $\mu(x,y)=xy$.}
\end{defi}

In the sequel, for simplicity, we often write $xy$ instead of
$\mu(x,y)$ since confusion rarely occurs. Note that Novikov algebras
are a special class of left-symmetric algebras with only equation
(\ref{Nov1}) satisfied. Define an associator on $\mathcal{A}$ by
\begin{eqnarray*}
a(x,y,z)=(xy)z-x(yz),\hspace{0.3cm}\mbox{for all}\
x,y,z\in\mathcal{A}.
\end{eqnarray*}
Then equation (\ref{Nov1}) is equivalent to
\begin{eqnarray*}
a(x,y,z)=a(y,x,z),\hspace{0.3cm}\mbox{for all}\ x,y,z\in\mathcal{A}.
\end{eqnarray*}
In other words, the associator is symmetric in $x,y$, from which the
left-symmetric algebra takes its name.

The commutator of a Novikov algebra (or a left-symmetric algebra)
$\mathcal {A}$
$$[x,y]=xy-yx,\hspace{0.3cm}\mbox{for all} \ x,y\in\mathcal{A},$$
defines a Lie algebra $Lie(\mathcal {A})$, which is called the {\it
sub-adjacent Lie algebra} of $\mathcal {A}$. Hence both Novikov
algebras and left-symmetric algebras are Lie-admissible algebras.

In order to establish tensor theory of Novikov algebras, Xu in
\cite{Xu2} introduced Novikov-Poisson algebras, which play important
roles in the study of simple Novikov algebras.
\begin{defi} {\rm A Novikov-Poisson algebra is a vector space $\mathcal {A}$ with two operations  $\cdot$ and
$\ast$ such that $(\mathcal{A},\cdot)$ forms a commutative
associative algebra and $(\mathcal{A},\ast)$ forms a Novikov algebra
for which
\begin{eqnarray}
(x\cdot y)\ast z&=&x\cdot(y\ast z),\label{NP-1}\\
(x\ast y)\cdot z-x\ast (y\cdot z)&=&(y\ast x)\cdot z-y\ast (x\cdot
z),\label{NP-2}
\end{eqnarray}
for all $x,y,z\in\mathcal{A}$.}
\end{defi}

 For an endomorphism of a Novikov-Poisson algebra $(\mathcal{A},\cdot,\ast)$, we mean a linear map
$\alpha:\mathcal{A}\longrightarrow\mathcal{A}$ preserving the two
operations, i.e.,
\begin{eqnarray*}
\alpha(x\cdot y)=\alpha(x)\cdot\alpha(y),\hspace{0.3cm} \alpha(x\ast
y)=\alpha(x)\ast\alpha(y), \hspace{0.3cm}\mbox{for all}\
x,y\in\mathcal {A}.
\end{eqnarray*}

The notion of Hom-Lie algebras was initially introduced by Hartwig,
Larsson and Silvestrov in\cite{HLS} motivated by examples of
deformed Lie algebras coming from twisted discretizations of vector
fields. A slightly more general definition of Hom-Lie algebras was
given by Makhlouf and Silvestrov in \cite{MS}, where the
Hom-associative algebras were also introduced.

\begin{defi}{\rm
 A Hom-Lie algebra is a triple
$(L,[\cdot,\cdot],\alpha)$ consisting of a vector space $L$, a
bilinear  map  $[ \cdot,\cdot]:L\times L\longrightarrow L$ and a
linear map $\alpha:L\longrightarrow L$ satisfying the following two
conditions
\begin{eqnarray*}
&&[x,y]=-[y,x],\hspace{0.5cm}\mbox{(skew-symmetry)}\\
&&[[x,y],\alpha(z)]+[[z,x],\alpha(y)]+[[y,z],\alpha(x)]=0,\hspace{0.5cm}\mbox{(Hom-Jacobi
\  identity)}
\end{eqnarray*}
 for all
$x,y,z\in L.$ }\end{defi}

Clearly, Lie algebras are spacial classes of Hom-Lie algebras when
$\alpha$ is the identity map. If the linear map $\alpha$ is an
algebra homomorphism with respect to the bracket $[\cdot,\cdot]$,
then the Hom-Lie algebras are a spacial case of quasi-hom-Lie
algebras \cite{LS}. Moreover, all these classes are special cases of
the more general quasi-Lie algebras \cite{LS2}, in which Lie color
algebras and in particular Lie superalgebras are also included.

\begin{defi} {\rm A  Hom-associative algebra is a triple $(V,
\mu, \alpha)$ consisting of a linear space $V$, a bilinear map
$¦Ì\mu : V \times V \longrightarrow V$ and a linear map $\alpha: V
\longrightarrow V$ satisfying
$$ \alpha(x)(yz) = (xy)\alpha(z),\hspace{0.3cm}\mbox{for all}\  x,y,z\in\mathcal {A}.$$ }
\end{defi}

We say that a Hom-associative algebra $(\mathcal {A},\mu,\alpha)$ is
{\it commutative }if $xy=yx$ holds for all $x,y\in\mathcal{A}$. The
following result provides a construction of Hom-Lie algebras from
Hom-associative algebras, extending the fundamental construction of
Lie algebras from associative algebras via commutator bracket
multiplication.
\begin{prop} {\rm(see [\ref{MS}, Proposition 1.6])} To any Hom-associative algebra defined by
a multiplication $\mu$ and a linear map $\alpha$ over an
$\mathbb{F}$-linear space $V$, one may associate a Hom-Lie algebra
defined by the bracket
\begin{eqnarray*}
[x,y]=xy-yx,\hspace{0.3cm}\mbox{for all}\ x,y\in V.
\end{eqnarray*}
\end{prop}

Following the patterns of Hom-Lie algebras and Hom-associative
algebras, Yau in\cite{Yau} introduced a twisted generalization of
Novikov algebras, called Hom-Novikov algebras.
\begin{defi}{\rm A
Hom-Novikov algebra is a triple $(\mathcal {A}, \mu, \alpha)$
consisting of  a vector space $\mathcal {A}$, a bilinear map $¦Ì\mu
: \mathcal {A} \times \mathcal {A} \longrightarrow \mathcal {A} $
and an algebra homomorphism $\alpha: \mathcal {A} \longrightarrow
\mathcal {A}$ satisfying
\begin{eqnarray}
(xy)\alpha(z) &=&(xz)\alpha(y),\label{Hom-Nov2}\\
(xy)\alpha(z)-\alpha(x)(yz) &=& (yx)\alpha(z)
-\alpha(y)(xz),\label{Hom-Nov1}
\end{eqnarray} for all $x,y,z\in\mathcal{A}.$
}\end{defi}

Clearly, Novikov algebras are examples of Hom-Novikov algebras by
setting $\alpha=\rm{id}$, where $\rm{id}$ is the identity map. It
was shown in \cite{Yau} that any Novikov algebra can be deformed
into a Hom-Novikov algebra along with an algebra endomorphism.

\begin{prop}\label{prop11}{\rm (see  [\ref{Yau}, Theorem 1.1]\,)}~Let $(\mathcal
{A}, \mu)$ be a Novikov algebra and $\alpha:\mathcal
{A}\longrightarrow \mathcal {A}$ be an algebra homomorphism. Then
$(\mathcal {A},\alpha\circ \mu,\alpha)$ is a Hom-Novikov algebra,
where $\circ$ denotes composition of functions.
\end{prop}

Similar to Hom-Novikov algebras, Yau in \cite{Yau5} introduced
Hom-Novikov-Poisson algebras.
\begin{defi}\label{Hom-NP}{\rm A Hom-Novikov-Poisson algebra is a
quadruple $(\mathcal {A},\cdot,\ast,\alpha)$, where $\mathcal {A}$
is a vector space, $\cdot$ and $\ast$ are bilinear maps and $\alpha$
is an algebra homomorphism of $\mathcal{A}$, such that $(\mathcal
{A},\ast,\alpha)$ forms a Hom-Novikov algebra, $(\mathcal
{A},\cdot,\alpha)$ forms a commutative Hom-associative algebra and
the following two compatible conditions hold:
\begin{eqnarray}
(x\cdot y)\ast\alpha(z)&=&\alpha (x)\cdot (y\ast z),\label{Hom-NP-1}\\
(x\ast y)\cdot\alpha(z)-\alpha(x)\ast(y\cdot z)&=&(y\ast
x)\cdot\alpha(z)-\alpha(y)\ast(x\cdot z),\label{Hom-NP-2}
\end{eqnarray}
for all $x,y,z\in\mathcal{A}$.}
\end{defi}

We can recover Novikov-Poisson algebras from Hom-Novikov-Poisson
algebras when $\alpha={\rm id}$. Some more interesting examples of
Hom-Novikov-Poisson algebras will be constructed in Section 5.

\vs{8pt}\

\cl{\bf\S3. \ Hom-Novikov
algebras}\setcounter{section}{3}\setcounter{equation}{0}\setcounter{theo}{0}
\vs{8pt} In this section, we establish the relationships between
Hom-Novikov algebras and Novikov algebras as well as the relations
between Hom-Novikov algebras and Hom-Lie algebras. Also, we provide
a construction of Hom-Novikov algebras from a commutative
Hom-associative algebra along with a suitable linear map but not a
derivation. Thus, it is different from that given in \cite{Yau}.

Just as a Novikov algebra can be deformed into a Lie algebra via the
commutator bracket, we have the similar result for Hom-Novikov
algebras and Hom-Lie algebras,  because in Hom-Novikov algebras hold
Jacobi-like identities.
\begin{prop}
Let $(\mathcal {A}, \mu, \alpha)$  be a Hom-Novikov algebra. For all
$x,y,z\in\mathcal {A}$, one has
\begin{eqnarray}
[x,y]\alpha(z)+[y,z]\alpha(x)+[z,x]\alpha(y)&=&0,\label{J1}\\
\alpha(x)[y,z]+\alpha(y)[z,x]+\alpha(z)[x,y]&=&0,\label{J2}
\end{eqnarray}
where $\mu(x,y)=xy$ and $[x,y]=xy-yx$.
\end{prop}
\noindent{\it Proof.~}For all $x,y,z\in\mathcal{A}$, we have
\begin{eqnarray*}
[x,y]\alpha(z)&=&(xy)\alpha(z)-(yx)\alpha(z),\\
{[y,z]}\alpha(x)&=&(yz)\alpha(x)-(zy)\alpha(x),\\
{[z,x]}\alpha(y)&=&(zx)\alpha(y)-(xz)\alpha(y).
\end{eqnarray*}
Since equation (\ref{Hom-Nov2}) holds for all $x,y,z\in\mathcal
{A}$, we have
\begin{eqnarray*}
&&[x,y]\alpha(z)+[y,z]\alpha(x)+[z,x]\alpha(y)\\
&=&\big((xy)\alpha(z)-(xz)\alpha(y)\big)+\big((yz)\alpha(x)-(yx)\alpha(z)\big)+\big((zx)\alpha(y)-(zy)\alpha(x)\big)\\
&=&0,
\end{eqnarray*}
which proves equation (\ref{J1}). Now using equations
(\ref{Hom-Nov1}) and (\ref{J1}), we have
\begin{eqnarray*}
&&\alpha(x)[y,z]+\alpha(y)[z,x]+\alpha(z)[x,y]\\
&=&\alpha(x)(yz)-\alpha(x)(zy)+\alpha(y)(zx)-\alpha(y)(xz)+\alpha(z)(xy)-\alpha(z)(yx)\\
&=&\big(\alpha(x)(yz)-\alpha(y)(xz)\big)+\big(\alpha(y)(zx)-\alpha(z)(yx)\big)+\big(\alpha(z)(xy)-\alpha(x)(zy)\big)\\
&=&\big((xy)\alpha(z)-(yx)\alpha(z)\big)+\big((yz)\alpha(x)-(zy)\alpha(x)\big)+\big((zx)\alpha(y)-(xz)\alpha(y)\big)\\
&=&[x,y]\alpha(z)+[y,z]\alpha(x)+[z,x]\alpha(y)\\
&=&0,
\end{eqnarray*}
which proves equation (\ref{J2}) and the proposition. \QED

\begin{coro} \label{c1} Let $(\mathcal
{A}, \mu, \alpha)$  be a Hom-Novikov algebra. Define a bilinear map
$[\cdot,\cdot]:\mathcal {A}\times\mathcal {A}\longrightarrow\mathcal
{A}$ by
$$[x,y]=\mu(x,y)-\mu(y,x),\hspace{0.3cm} \mbox{for all}\ x,y\in\mathcal {A}.$$
Then $HLie(\mathcal {A})=(\mathcal {A}, [\cdot,\cdot], \alpha)$ is a
Hom-Lie algebra with the same twisting map $\alpha$.
\end{coro}

The Lie admissible algebras were introduced by A.A. Albert in 1948.
Makhlouf and Silvestrov in \cite{MS} extended the notions and
results about Lie admissible algebras to Hom-algebras, called
Hom-Lie admissible algebras. According to Corollary \ref{c1}, any
Hom-Novikov algeba is Hom-Lie admissible with the same twisting map.
We call such a Hom-Lie algebra $HLie(\mathcal {A})$ the {\it
sub-adjacent Hom-Lie algebra }of the Hom-Novikov algebra $(\mathcal
{A}, \mu, \alpha)$.

\begin{defi} {\rm Let $(\mathcal
{A}, \mu, \alpha)$  be a Hom-Novikov algebra, which is called
\begin{itemize}\item[\rm(i)]{\it regular} if $\alpha$ is an algebra
automorphism;
\item[\rm(ii)]{\it involutive} if $\alpha$ is an involution, i.e.,
$\alpha^2=\rm{id}.$
\end{itemize}}
\end{defi}

Yau in \cite{Yau} gave a way to construct Hom-Novikov algebras,
starting from a Novikov algebra and an algebra endomorphism. In the
following, we provide a construction of Novikov algebras from
Hom-Novikov algebras along with an algebra automorphism.
\begin{prop} \label{prop22}If $(\mathcal
{A}, \mu, \alpha)$  is an involutive Hom-Novikov algebra, then
$(\mathcal {A}, \alpha\circ\mu)$ is a Novikov algebra.
\end{prop}
\noindent{\it Proof.~} For convenience, we write $\mu(x,y)=xy$ and
$x\star y=\alpha(xy),$ for all $x,y\in\mathcal {A}$. Hence, it needs
to show
\begin{eqnarray}
(x\star y)\star z &=&(x\star z)\star y,\label{Nov4}\\
(x\star y)\star z-x\star (y\star z) &=& (y\star x)\star z -y\star
(x\star z),\label{Nov3}
\end{eqnarray}
for all $x,y,z\in\mathcal {A}$. Since $\alpha$ is an involution, we
have
\begin{eqnarray*}
(x\star y)\star z&=&\alpha\big(\alpha(xy)z\big)
=\big(\alpha^{2}(xy)\big)\alpha(z) =(xy)\alpha(z).
\end{eqnarray*}
Similarly, we have $(x\star z)\star y=(xz)\alpha(y)$, from which
equation (\ref{Nov4}) follows since $(\mathcal {A}, \mu, \alpha)$ is
a Hom-Novikov algebra. Furthermore, using equation (\ref{Hom-Nov1}),
we have
\begin{eqnarray*}
(x\star y)\star z-x\star (y\star
z)&=&\alpha\big(\alpha(xy)z\big)-\alpha\big(x\alpha(yz)\big)\\
&=&(xy)\alpha(z)-\alpha(x)(yz)\\
&=&(yx)\alpha(z)-\alpha(y)(xz)\\
&=&\alpha\big(\alpha(yx)z\big)-\alpha\big(y\alpha(xz)\big)\\
&=&(y\star x)\star z -y\star (x\star z),
\end{eqnarray*}
which proves equation (\ref{Nov3}) and the proposition.\QED
\begin{prop} If $(\mathcal
{A}, \mu, \alpha)$  is a regular Hom-Novikov algebra,
 then $(\mathcal {A},[\cdot,\cdot]_{\alpha^{-1}}=
\alpha^{-1}\circ [\cdot,\cdot])$ is a Lie algebra, where
$[x,y]=\mu(x,y)-\mu(y,x)=xy-yx,$ for all $x,y\in\mathcal {A}.$ In
particular, if $\alpha$ is an involution, then  $(\mathcal {A},
[\cdot,\cdot]_{\alpha}=\alpha\circ [\cdot,\cdot])$ is a Lie algebra.
\end{prop}
\noindent{\it Proof.~} For any $x,y,z\in\mathcal{A}$, we have
\begin{eqnarray*}
[[x,y]_{\alpha^{-1}},z]_{\alpha^{-1}}
&=&[\alpha^{-1}(xy-yx),z]_{\alpha^{-1}}\\
&=&\alpha^{-1}\big(\alpha^{-1}(xy-yx)z-z\alpha^{-1}(xy-yx)\big)\\
&=&\alpha^{-2}\big((xy-yx)\alpha(z)-\alpha(z)(xy-yx)\big)\\
&=&\alpha^{-2}\big([x,y]\alpha(z)-\alpha(z)[x,y]\big).
\end{eqnarray*}
Similarly, we have
\begin{eqnarray*}
[[y,z]_{\alpha^{-1}},x]_{\alpha^{-1}}&=&
\alpha^{-2}\big([y,z]\alpha(x)-\alpha(x)[y,z]\big),\\
{[[z,x]_{\alpha^{-1}},y]_{\alpha^{-1}}}&=&
\alpha^{-2}\big([z,x]\alpha(y)-\alpha(y)[z,x]\big).\\
\end{eqnarray*}
Then it follows from equations (\ref{J1}) and  (\ref{J2}) that
$$[[x,y]_{\alpha^{-1}},z]_{\alpha^{-1}}+[[y,z]_{\alpha^{-1}},x]_{\alpha^{-1}}
+{[[z,x]_{\alpha^{-1}},y]_{\alpha^{-1}}}=0.$$ Clearly,
$[x,y]=-[y,x].$ Since $\alpha^{-1}$ is an automorphism, we have
$$[x,y]_{\alpha^{-1}}=-[y,x]_{\alpha^{-1}},$$
which proves that $(\mathcal {A},[\cdot,\cdot]_{\alpha^{-1}})$ is a
Lie algebra. It follows immediately that $(\mathcal
{A},[\cdot,\cdot]_{\alpha})$ is also a Lie algebra since
$\alpha=\alpha^{-1}$ when $\alpha$ is an involution.\QED\vskip7pt

Let $(\mathcal {A},\mu)$ be a commutative associative algebra and
$D$ be a derivation of $(\mathcal {A},\mu)$. Then the new product
\begin{eqnarray}\label{Gd1}
x\star_\lambda y=xD(y)+\lambda x y, \hspace{0.3cm}\mbox{for all}\
x,y\in\mathcal {A},
\end{eqnarray}
makes $(\mathcal {A},\star_\lambda)$ become a Novikov algebra for
$\lambda=0$ by Gel'fand and Dorfman \cite{GD}, for
$\lambda\in\mathbb{F}$ by Filipov \cite{F}, and for a fixed element
$\lambda\in\mathcal {A}$ by Xu \cite{Xu1}. Yau \cite{Yau}
generalized this construction to Hom-Novikov algebras $(\mathcal
{A},\mu, \alpha)$ in $\lambda=0$ setting together with an additional
condition that $D$ commutes with $\alpha$ . It is natural to
consider other type of linear maps to replace the derivation $D$ in
(\ref{Gd1}), such that  $(\mathcal{A},\star_\lambda)$ forms a
Novikov algebra or a Hom-Novikov algebra.

Let $(\mathcal {A},\mu,\alpha)$ be a commutative Hom-associative
algebra with a linear selfmap $\partial$ commuting with $\alpha$.
Consider the following operation on $\mathcal{A}$:
\begin{eqnarray}\label{Gd-alpha1}
x\star y=x\partial(y), \hspace{0.3cm} \mbox{for all}\ x,y\in\mathcal
{A}.
\end{eqnarray}
For all $x,y,z\in \mathcal {A}$, by straightforward calculation, we
have
\begin{eqnarray*}
(x\star
y)\star\alpha(z)&=&(x\partial(y))\star\alpha(z)\\
&=&(x\partial(y))\partial(\alpha(z))\\
&=&(x\partial(y))\alpha(\partial(z))\\
&=&\alpha(x)(\partial(y)\partial(z)),\\
\end{eqnarray*}
and
\begin{eqnarray*}
(x\star
z)\star\alpha(y)&=&(x\partial(z))\star\alpha(y)\\
&=&(x\partial(z))\partial(\alpha(y))\\
&=&(x\partial(z))\alpha(\partial(y))\\
&=&\alpha(x)(\partial(z)\partial(y))\\
&=&\alpha(x)(\partial(y)\partial(z)),\\
\end{eqnarray*}
where both commutativity and Hom-associativity of $(\mathcal
{A},\mu,\alpha)$ are used. Hence, we have
\begin{eqnarray*}
(x\star y)\star\alpha(z)=(x\star z)\star\alpha(y).
\end{eqnarray*}
Similarly, we have
\begin{eqnarray}\label{hn1}
(x\star y)\star\alpha(z)-\alpha(x)\star(y\star z)
=\alpha(x)(\partial(y)\partial(z))-\alpha(x)\partial(y\partial(z)),
\end{eqnarray}
and
\begin{eqnarray}\label{hn2}
(y\star x)\star\alpha(z)-\alpha(y)\star(x\star
z)=\alpha(y)(\partial(x)\partial(z))-\alpha(y)\partial(x\partial(z)).
\end{eqnarray}
If, in addition, $\partial$ satisfies the following condition
\begin{eqnarray}\label{Gd-2}
\partial(x\partial (y))=\partial(x)\partial (y), \hspace{0.3cm}\mbox{for all}\ x,y\in\mathcal{A}.
\end{eqnarray}
Then it follows from equations (\ref{hn1}) and (\ref{hn2}) that
\begin{eqnarray*}
(x\star y)\star\alpha(z)-\alpha(x)\star(y\star z)=(y\star
x)\star\alpha(z)-\alpha(y)\star(x\star z)=0.
\end{eqnarray*}
Now from the discussions above, we obtain
\begin{theo}\label{twist}
Let $(\mathcal {A},\mu, \alpha)$ be a commutative Hom-associative
algebra with a linear selfmap $\partial$ satisfying equation
(\ref{Gd-2}) and commuting with  $\alpha$. Then $(\mathcal
{A},\star, \alpha)$ is a Hom-Novikov algebra, where $\star$ is
defined by (\ref{Gd-alpha1}).
\end{theo}

We have the following corollary by setting $\alpha={\rm id}$ in the
theorem.
\begin{coro}\label{c2} Let $(\mathcal {A},\mu)$ be a commutative associative
algebra with a linear selfmap $\partial$ satisfying condition
(\ref{Gd-2}). Then $(\mathcal {A},\star)$ is a Novikov algebra in
which $\star$ is defined by (\ref{Gd-alpha1}).
\end{coro}

\begin{exam}{\rm Let  $(\mathcal{A}=\mathcal{A}_{0}\oplus \mathcal{A}_{
1},\mu)$ be a complex superalgebra, where $\mathcal{A}_{
0}=\mathbb{C}[t,t^{-1}]$ is the Laurent polynomials in one variable
and $\mathcal{A}_{1}=\theta\mathbb{C}[t,t^{-1}]$ in which $\theta$
is the Grassman variable satisfying $\theta^2=0$. We assume that $t$
commutes with $\theta$. The generators of $\mathcal{A}$ are of the
form $t^n$ and $\theta t^n$ for all $n\in \mathbb{Z}$. Then
$\mathcal{A}$ is a commutative associative  superalgebra under the
usual multiplication.

Define a linear map $\partial:\mathcal{A}\longrightarrow\mathcal{A}$
by
\begin{eqnarray}\label{d}
\partial(t^n+\theta t^m)=t^n,\hspace{0.3cm}\mbox{for all}\ m,n\in
\mathbb{Z}.
\end{eqnarray}
 It is easy to check that $\partial$ satisfies condition
(\ref{Gd-2}).  According to Corollary \ref{c2},
$(\mathcal{A},\star_1)$ is a Novikov algebra, where the operation
$\star_1$ is defined by
\begin{eqnarray*}
(t^{n_1}+\theta t^{m_1} )\star_1(t^{n_2}+\theta t^{m_2}
)=(t^{n_1}+\theta t^{m_1})\partial (t^{n_2}+\theta
t^{m_2})=t^{n_1+n_2}+\theta t^{m_1+n_2},
\end{eqnarray*}
for all $n_1,n_2,m_1,m_2\in\mathbb{Z}$.

Furthermore, define an algebra endomorphism $\alpha$ on
$\mathcal{A}$ by
$$\alpha(t^n+\theta t^m)=(t+c)^n,\hspace{0.3cm}\mbox{for all}\  m,n\in
\mathbb{Z},$$ where $c$ is a fixed element in $\mathbb{C}$. Then
$(\mathcal{A},\alpha\circ\mu,\alpha)$ is a commutative
Hom-associative algebra.  For convenience, we write $\cdot$ instead
of $\alpha\circ\mu$, then we have
\begin{eqnarray*}
(t^{n_1}+\theta t^{m_1} )\cdot (t^{n_2}+\theta t^{m_2}
)=(t+c)^{n_1+n_2},
\end{eqnarray*}
for all $n_1,n_2,m_1,m_2\in\mathbb{Z}$. Let $\partial$ be as defined
by (\ref{d}). Clearly, $\partial$ commutes with $\alpha$. Moreover,
we have
\begin{eqnarray*}
\partial\big((t^{n_1}+\theta t^{m_1} )\cdot \partial(t^{n_2}+\theta t^{m_2}
)\big)=(t+c)^{n_1+n_2},
\end{eqnarray*}
and
\begin{eqnarray*}
\partial(t^{n_1}+\theta t^{m_1} )\cdot \partial(t^{n_2}+\theta t^{m_2}
)=(t+c)^{n_1+n_2}.
\end{eqnarray*}
Hence, $\partial$ satisfies condition (\ref{Gd-2}). Define a new
bilinear operation on $\mathcal{A}$ by
\begin{eqnarray*}\label{d3}
(t^{n_1}+\theta t^{m_1} )\star_2(t^{n_2}+\theta t^{m_2}
)=(t^{n_1}+\theta t^{m_1})\cdot \partial (t^{n_2}+\theta
t^{m_2})=(t+c)^{n_1+n_2},
\end{eqnarray*}
for all $n_1,n_2,m_1,m_2\in\mathbb{Z}$. Then
$(\mathcal{A},\star_2,\alpha)$ forms a Hom-Novikov algebra by
Theorem \ref{twist}. }\end{exam}

 \vskip7pt
 \cl{\bf\S4. \ Quadratic Hom-Novikov
Algebras}\setcounter{section}{4}\setcounter{equation}{0}\setcounter{theo}{0}
\vskip8pt

In this section, we extend the notion of quadratic Novikov algebra
 to quadratic Hom-Novikov algebras and provide some properties.

Recall that a  quadratic Lie algebra is a triple $(\mathcal
{G},[\cdot,\cdot],B)$, where $(\mathcal {G},[\cdot,\cdot])$ is a Lie
algebra and $B:\mathcal {G}\times \mathcal
{G}\longrightarrow\mathbb{F}$ is a symmetric nondegenerate bilinear
form satisfying
\begin{eqnarray}\label{F1}
B([x,y],z)=B(x,[y,z]),\hspace{0.3cm}\mbox{for all}\ x,y,z\in\mathcal
{G},
\end{eqnarray}
which is called the {\it invariance} of $B$. It is easy to see that
identity (\ref{F1}) is equivalent to
\begin{eqnarray*}
B([x,y],z)=-B(y,[x,z]),\hspace{0.3cm}\mbox{for all}\
x,y,z\in\mathcal {G}.
\end{eqnarray*}
Assume that $(\mathcal {A},\mu)$ is a Novikov algebra. A bilinear
form $B:\mathcal {A}\times\mathcal {A}\longrightarrow \mathbb{F}$ is
said to be {\it invariant} or {\it associative} if and only if
$$B(xy,z)=B(x,yz),\hspace{0.3cm}\mbox{for all}\ x,y,z\in\mathcal {A}.$$
A Novikov algebra  $(\mathcal {A},\mu)$ with a symmetric
nondegenerate invariant  bilinear form is called a {\it quadratic
Novikov algebra} and denoted by  $(\mathcal {A},\mu,B)$.

Benayadi and Makhlouf in \cite{BM} extended the notion of quadratic
Lie algebra to Hom-Lie algebras and obtained quadratic Hom-Lie
algebras.

\begin{defi} {\rm A quadratic Hom-Lie algebra is a quadruple $(\mathcal {A},[\cdot,\cdot],\alpha,B)$
such that $(\mathcal {A},[\cdot,\cdot],\alpha)$ is a Hom-Lie algebra
with a symmetric invariant nondegenerate bilinear form $B$
satisfying
\begin{eqnarray}\label{F2}
B(\alpha(x),y)=B(x,\alpha(y)),\hspace{0.3cm}\mbox{for all}\
x,y\in\mathcal{A}.
\end{eqnarray}}
\end{defi}

We can define quadratic Hom-Novikov algebras as follows.
\begin{defi}{\rm A quadratic Hom-Novikov algebra $(\mathcal {A},\mu,\alpha,B)$
is a Hom-Novikov algebra $(\mathcal {A},\mu,\alpha)$  with a
 symmetric nondegenerate bilinear form satisfying
\begin{eqnarray}\label{quadratic}
B(xy,\alpha(z))=B(\alpha(x),yz),\hspace{0.3cm}\mbox{for all}\
x,y,z\in\mathcal {A}.
\end{eqnarray}}
\end{defi}

We recover quadratic Novikov algebras when $\alpha={\rm id}$. The
following result says the sub-adjacent Hom-Lie algebra of a
quadratic Hom-Novikov algebra is also quadratic.
\begin{prop} \label{p4}Let $(\mathcal {A},\mu,\alpha,B)$ be a quadratic Hom-Novikov
algebra and $HLie(\mathcal{A})=(\mathcal{A}, [\cdot,\cdot],\alpha)$
be the sub-adjacent Hom-Lie algebra of $\mathcal{A}$. If $\alpha$ is
an automorphism satisfying
\begin{eqnarray}\label{F3}
B(\alpha(x),y)=B(x,\alpha(y)), \hspace{0.3cm}\mbox{for all}\
x,y\in\mathcal{A}.
\end{eqnarray}
Then $(\mathcal{A}, [\cdot,\cdot],\alpha, B_\alpha)$ is a quadratic
Hom-Lie algebra, where $B_\alpha(x,y)=B( \alpha(x),y)$, for all
$x,y\in\mathcal{A}$.
\end{prop}
\noindent{\it Proof.~} Since $B$ is a nondegenerate bilinear form
and $\alpha$ is an automorphism, $B_\alpha$ is a nondegenerate
bilinear
 form on $\mathcal{A}$. For all $x,y,z\in\mathcal{A}$,
using the properties of $B$, we have
\begin{eqnarray*}
B_\alpha([x,y],z)&=&B\big(\alpha([x,y]),z\big)
\\&=&B([x,y],\alpha(z))\\
&=&B(xy,\alpha(z))-B(yx,\alpha(z))\\
&=&B(\alpha(x),yz)-B(\alpha(x),zy)\\
&=&B(\alpha(x),[y,z])\\&=&B_\alpha(x,[y,z]).
\end{eqnarray*}
Hence $B_\alpha$ is invariant. Using symmetry of $B$ and equation
(\ref{F3}), we have
\begin{eqnarray*}
B_\alpha(x,y)=B(\alpha(x),y)=B(y,\alpha(x))=B(\alpha(y),x)=B_\alpha(y,x),
\end{eqnarray*}
which proves $B_\alpha$ is symmetric. Using equation (\ref{F3})
again, we have
$$B_\alpha(\alpha(x),y)=B(\alpha(\alpha(x)),y)=B(\alpha(x),\alpha(y))=B_\alpha(x,\alpha(y)),$$
which completes the proof.\QED
\begin{coro} Let $(\mathcal{A},\mu, B)$ be a quadratic Novikov
algebra with an automorphism $\alpha$
 satisfying equation (\ref{F3}) and $(\mathcal{A},[\cdot,\cdot])$ be the
 sub-adjacent Lie algebra. Then $(\mathcal{A},
[\cdot,\cdot]_\alpha=\alpha\circ[\cdot,\cdot],\alpha, B_\alpha)$
forms a quadratic Hom-Lie algebra, where $B_\alpha(x,y)=B(
\alpha(x),y)$, for all $x,y\in\mathcal{A}$.
\end{coro}
\noindent{\it Proof.~} It follows from [\ref{Yau2}, Corollary 2.6]
that $(\mathcal{A}, [\cdot,\cdot]_\alpha,\alpha)$ is a Hom-Lie
algebra.  Using the similar arguments as those in the proof of
Proposition \ref{p4}, we get $B_{\alpha}$ is a symmetric
 nondegenerate bilinear form with equation (\ref{F2}) satisfied. It
remains to show that $B_{\alpha}$ is invariant. For all
$x,y,z\in\mathcal{A}$, using invariance and symmetry of $B$, we have
\begin{eqnarray*}
B_\alpha([x,y]_{\alpha},z)&=&B\big(\alpha([x,y]_{\alpha}),z\big)
\\&=&B\big([x,y]_{\alpha},\alpha(z)\big)\\
&=&B(\alpha(x)\alpha(y),\alpha(z))-B(\alpha(y)\alpha(x),\alpha(z))\\
&=&B(\alpha(x),\alpha(y)\alpha(z))-B(\alpha(z)\alpha(y),\alpha(x))\\
&=&B(\alpha(x),[y,z]_{\alpha})\\&=&B_\alpha(x,[y,z]_\alpha),
\end{eqnarray*}
which proves the invariance of $B_{\alpha}$ and the result.\QED
\begin{prop}Let $(\mathcal {A},\mu, \alpha, B)$ be a quadratic Hom-Novikov
algebra, where $\alpha$ is an involution satisfying equation
(\ref{F3}). Then $(\mathcal {A},\star=\alpha\circ\mu,B)$ is  a
quadratic Novikov algebra.
\end{prop}
\noindent{\it Proof.~} $({A},\star=\alpha\circ\mu,\alpha)$ is a
Hom-Novikov algebra by Proposition \ref{prop22}. It suffices to show
that $B$ is invariant under the operation $\star$. For all
$x,y,z\in\mathcal{A}$, we have
\begin{eqnarray*}
B(x,y\star
z)=B(x,\alpha(y)\alpha(z))=B(\alpha(x),yz)=B(xy,\alpha(z))=B(\alpha(x)\alpha(y),z)=B(x\star
y,z),
\end{eqnarray*}
which completes the proof.\QED

\begin{prop}Let $(\mathcal {A},\mu,\alpha, B)$ be a quadratic Hom-Novikov
algebra, where $\alpha$ is an automorphism satisfying equation
(\ref{F3}). Then, $(\mathcal{A},\star=\alpha\circ\mu,\alpha^{2},
B_{\alpha^2})$ is a quadratic Hom-Novikov algebra, where
$B_{\alpha^2}(x,y)=B(\alpha^2(x),y)$, for all $x,y\in\mathcal{A}$.
\end{prop}
\noindent{\it Proof.~}It follows from [\ref{Yau5}, Corollary 2.12])
that $(\mathcal{A},\alpha\circ\mu,\alpha^{2})$ forms a Hom-Novikov
algebra. Since $B$ is a nondegenerate bilinear form on $\mathcal{A}$
and $\alpha$ is an automorphism,  $B_{\alpha^2}$ is a nondegenerate
bilinear form. For all $x,y,z\in\mathcal{A}$, using the hypothesis,
we have
$$B_{\alpha^2}(x, y)=B(\alpha^2(x),y)=B(\alpha(x),\alpha(y))=B(x,\alpha^2(y))=B(\alpha^2(y),x)=B_{\alpha^2}(y,x).$$
Thus, $B_{\alpha^2}$ is symmetric. Moreover, we have
\begin{eqnarray*}
B_{\alpha^2}\big(\alpha^{2}(x),y\star z\big)&=&B\big(\alpha^{4}(x),\alpha(y)\alpha(z)\big)\\
&=&B\big(\alpha^{3}(x),\alpha^2(y)\alpha^2(z)\big)\\
&=&B\big(\alpha^{2}(x)\alpha^2(y),\alpha^3(z)\big)\\
&=&B\big(\alpha(x)\alpha(y),\alpha^4(z)\big)\\
&=&B_{\alpha^2}\big(x\star y,\alpha^2(z)\big),
\end{eqnarray*}
which proves the invariance of $B_{\alpha^2}$ and the
proposition.\QED
\begin{coro}Let $(\mathcal {A},\mu,\alpha, B)$ be a quadratic Hom-Novikov
algebra and $\alpha$ be an automorphism satisfying equation
(\ref{F3}). Then $(\mathcal{A},\alpha^n\circ\mu,\alpha^{n+1},
B_{\alpha^{n+1}})$ is a quadratic Hom-Novikov algebra, where
$B_{\alpha^{n+1}}(x,y)=B(\alpha^{n+1}(x),y)$, for all $n>0$ and
$x,y\in\mathcal{A}$.
\end{coro}

Let $(\mathcal{A},\mu,\alpha)$ be a Hom-Novikov algebra, whose
center is denoted by $\mathcal {Z}(\mathcal{A})$ and defined by
\begin{eqnarray*}
\mathcal {Z}(\mathcal{A})=\{x\in\mathcal{A}|xy=yx=0,\ \mbox{for
all}\ y\in\mathcal{A}\}.
\end{eqnarray*}
Let $(\mathcal {G},[\cdot,\cdot],\beta)$ be a Hom-Lie algebra. The
lower central series of $\mathcal {G}$ is defined as usual, i.e.,
$\mathcal {G}^1=\mathcal {G}$, $\mathcal {G}^{i+1}=[\mathcal
{G},\mathcal {G}^i]$ for all $i\geq 1$. We call $\mathcal {G}$ is
{\it $i$-step nilpotent }if $\mathcal {G}^{i+1}=0$. The center of
the Hom-Lie algebra is denoted by $\mathcal {C}(\mathcal {G})$ and
defined by
\begin{eqnarray*}
\mathcal {C}(\mathcal {G})=\{x\in \mathcal {G}|[x,y]=0,\ \mbox{for
all}\ y\in \mathcal {G}\}.
\end{eqnarray*}
\begin{prop}\label{prop4}Let $(\mathcal{A},\mu, \alpha,B)$ be a quadratic Hom-Novikov algebra and
$HLie(\mathcal{A})$ be the sub-adjacent Hom-Lie algebra. If $\alpha$
is an automorphism, then
$[HLie(\mathcal{A}),HLie(\mathcal{A})]\subseteq \mathcal
{Z}(\mathcal{A}).$ As a consequence, $HLie(\mathcal{A})$ is $2$-step
nilpotent.
\end{prop}
\noindent{\it Proof.~} For all $x,y,z,w\in\mathcal{A}$, using
equations (\ref{Hom-Nov2}), (\ref{Hom-Nov1}) and (\ref{quadratic}),
we have
\begin{eqnarray*}
B(\alpha(x)[y,z],\alpha^2(w))&=&B(\alpha^2(x),[y,z]\alpha(w))\\
&=&B(\alpha^2(x),(yz)\alpha(w)-(zy)\alpha(w))\\
&=&B(\alpha^2(x),\alpha(y)(zw)-\alpha(z)(yw))\\
&=&B(\alpha(x)\alpha(y),\alpha(zw))-B(\alpha(x)\alpha(z),\alpha(yw))\\
&=&B((xy)\alpha(z)-(xz)\alpha(y),\alpha^2(w))\\
&=&0.
\end{eqnarray*}
Using the symmetry of $B$ and equation (\ref{quadratic}), we have
$$B([y,z]\alpha(w),\alpha^2(x))=B(\alpha(x)[y,z],\alpha^2(w))=0,$$
 which implies
$[y,z]\subseteq \mathcal {Z}(\mathcal{A})$ since $\alpha$ is an
automorphism and $B$ is nondegenerate. Hence, we have
$[HLie(\mathcal{A}),HLie(\mathcal{A})]\subseteq \mathcal
{Z}(\mathcal{A}).$ Obviously, $\mathcal {Z}(\mathcal{A})\subseteq
\mathcal {C}(HLie(\mathcal{A}))$. Then it follows that
$HLie(\mathcal{A})$ is $2$-step nilpotent.\QED \vskip8pt

 \cl{\bf\S5. \ Hom-Novikov-Poisson
Algebras}\setcounter{section}{5}\setcounter{equation}{0}\setcounter{theo}{0}
\vskip8pt

In this section, we discuss Hom-Novikov-Poisson algebras. We give
different ways to  construct Hom-Novikov-Poisson algebras and
provide some interesting examples.

 The following result (see also [\ref{Yau5}, Corollary 2.15])says
 that any Hom-Novikov-Poisson algebra can be constructed from
 a Novikov-Poisson algebra and an algebra endomorphism.
\begin{prop}\label{prop3} Let $\big(\mathcal {A},\mu,\nu\big)$ be
a Novikov-Poisson algebra such that $\big(\mathcal {A},\nu\big)$ is
a commutative associative algebra and $\big(\mathcal {A},\mu\big)$
is a Novikov algebra. If $\alpha:\mathcal {A}\longrightarrow\mathcal
{A}$ is an algebra homomorphism. Then $\mathcal
{A}_{\alpha}=\big(\mathcal
{A},\mu_{\alpha}=\alpha\circ\mu,\nu_{\alpha}=\alpha\circ\nu,\alpha\big)$
forms a Hom-Novikov-Poisson algebra.
\end{prop}
\noindent{\it Proof.~}It follows that
$\big(\mathcal{A},\mu_{\alpha},\alpha\big)$ is a Hom-Novikov algebra
from Proposition \ref{prop11} and
$\big(\mathcal{A},\nu_{\alpha},\alpha\big)$ is a Hom-associative
algebra from [\ref{Yau2}, Corollary 2.6]. Clearly,
$\big(\mathcal{A},\nu_{\alpha},\alpha\big)$ is commutative
 since $\alpha$ is an algebra homomorphism and $(\mathcal{A},\nu)$ is
 commutative.  By Definition \ref{Hom-NP}, it suffices to check conditions
(\ref{Hom-NP-1}) and (\ref{Hom-NP-2}). For all $x,y,z\in\mathcal
{A}$,  Using equation (\ref{NP-2}), we have
\begin{eqnarray*}
\nu_\alpha\big(\mu_\alpha(x,y),\alpha(z)\big)-\mu_\alpha\big(\alpha(x),\nu_\alpha(y,z)\big)&=&
\alpha\circ\nu\big(\alpha\circ\mu(x,y),\alpha(z)\big)-\alpha\circ\mu\big(\alpha(x),\alpha\circ\nu(y,z)\big)\\
&=&\alpha\circ\alpha\circ\nu\big(\mu(x,y),z\big)-\alpha\circ\alpha\circ\mu\big(x,\nu(y,z)\big)\\
&=&\alpha\circ\alpha\circ\Big(\nu\big(\mu(x,y),z\big)-\mu\big(x,\nu(y,z)\big)\Big)\\
&=&\alpha\circ\alpha\circ\Big(\nu\big(\mu(y,x),z\big)-\mu\big(y,\nu(x,z)\big)\Big)\\
&=&\alpha\circ\nu\big(\alpha\circ\mu(y,x),\alpha(z)\big)-\alpha\circ\mu\big(\alpha(y),\alpha\circ\nu(x,z)\big)\\
&=&\nu_\alpha\big(\mu_\alpha(y,x),\alpha(z)\big)-\mu_\alpha\big(\alpha(y),\nu_\alpha(x,z)\big),
\end{eqnarray*}
which proves equation (\ref{Hom-NP-2}). Using equation (\ref{NP-1}),
we have
\begin{eqnarray*}
\mu_\alpha\big(\nu_\alpha(x,y),\alpha(z)\big)&=&
\alpha\circ\mu\big(\alpha\circ\nu(x,y),\alpha(z)\big)\\
&=&\alpha\circ\alpha\circ\mu\big(\nu(x,y),z\big)\\
&=&\alpha\circ\alpha\circ\nu\big(x,\mu(y,z)\big)\\
&=&\nu_\alpha\big(\alpha(x),\mu_\alpha(y,z)\big),
\end{eqnarray*}
which proves equation (\ref{Hom-NP-1}) and the result.\QED\vskip7pt

Suppose that  $\big(\mathcal {A}_1,\mu^1,\nu^1\big)$ and
$\big(\mathcal {A}_2,\mu^2,\nu^2\big)$ are Novikov-Poisson algebras.
Define two operations $\nu^1\otimes\nu^2$ and $\mu^1\otimes\mu^2$ on
$\mathcal {A}_1\otimes\mathcal {A}_2$ by
\begin{eqnarray}
(\nu^1\otimes\nu^2)(x_1\otimes x_2,y_1\otimes
y_2)&=&\nu^1(x_1,y_1)\otimes\nu^2(x_2,y_2),\label{tensor1}\\
(\mu^1\otimes\mu^2)(x_1\otimes x_2,y_1\otimes
y_2)&=&\mu^1(x_1,y_1)\otimes\nu^2(x_2,y_2)+\nu^1(x_1,y_1)\otimes\mu^2(x_2,y_2),\label{tensor2}
\end{eqnarray}
for all $x_1, \ y_1\in\mathcal {A}_1$ and $x_2, \ y_2\in\mathcal
{A}_2$. Then $\big(\mathcal {A}_1\otimes \mathcal {A}_2,
\mu^1\otimes\mu^2,\nu^1\otimes\nu^2\big)$ forms a Novikov-Poisson
algebra (cf.\cite{Xu1, Xu2}). Additionally, assume that $\alpha$ and
$\beta$ are algebra homomorphisms of $\big(\mathcal
{A}_1,\mu^1,\nu^1\big)$ and $\big(\mathcal {A}_2,\mu^2,\nu^2\big)$,
respectively. It follows from Proposition \ref{prop3} that
$\big(\mathcal {A}_1,\mu^1_\alpha,\nu^1_\alpha,\alpha\big)$ and
$\big(\mathcal {A}_2,\mu^2_\beta,\nu^2_\beta,\beta\big)$ are
Hom-Novikov-Poisson algebras. Define a linear map
$\alpha\otimes\beta$ on $\mathcal {A}_1\otimes \mathcal {A}_2$ by
\begin{eqnarray}
(\alpha\otimes\beta)(x\otimes
y)=\alpha(x)\otimes\beta(y),\hspace{0.3cm} \mbox{for all} \
x\in\mathcal{A}_1,\ y\in\mathcal {A}_2.
\end{eqnarray}
 For all $x_1, y_1\in \mathcal
{A}_1$ and $x_2, y_2\in \mathcal {A}_2$, we have
\begin{eqnarray*}
(\alpha\otimes\beta)\big(\nu^1\otimes\nu^2(x_1\otimes x_2,
y_1\otimes y_2)\big)&=&(\alpha\otimes\beta)\big(\nu^1(x_1,y_1)\otimes\nu^2(x_2,y_2)\big)\\
&=&\alpha\circ\nu^1(x_1,y_1)\otimes\beta\circ\nu^2(x_2,y_2)\\
&=&\nu^1(\alpha (x_1),\alpha (y_1))\otimes\nu^2(\beta(x_2),\beta(y_2))\\
&=&(\nu^1\otimes\nu^2)\big(\alpha(x_1)\otimes\beta(x_2),\alpha(y_1)\otimes\beta(y_2)\big)\\
&=&(\nu^1\otimes\nu^2)\big(\alpha\otimes\beta(x_1\otimes
x_2),\alpha\otimes\beta(y_1\otimes y_2)\big).
\end{eqnarray*}
Similarly, we have
\begin{eqnarray*}
(\alpha\otimes\beta)\big(\mu^1\otimes\mu^2(x_1\otimes x_2,y_1\otimes
y_2)\big) =(\mu^1\otimes\mu^2)\big(\alpha\otimes\beta(x_1\otimes
x_2),\alpha\otimes\beta(y_1\otimes y_2)\big).
\end{eqnarray*}
Thus $\alpha\otimes\beta$ is an algebra homomorphism of
$\big(\mathcal {A}_1\otimes \mathcal {A}_2,\mu^1\otimes\mu^2,
\nu^1\otimes\nu^2\big)$.

On the other hand, consider the following two operations on
$\mathcal {A}_1\otimes \mathcal {A}_2$ :
\begin{eqnarray*}
(\nu^1_\alpha\otimes\nu^2_\beta)(x_1\otimes x_2,y_1\otimes
y_2)&=&\nu^1_\alpha(x_1,y_1)\otimes\nu^2_\beta(x_2, y_2),\\
(\mu^1_\alpha\otimes\mu^2_\beta)(x_1\otimes x_2,y_1\otimes
y_2)&=&\mu^1_\alpha(x_1,y_1)\otimes\nu^2_\beta(x_2,y_2)+\nu^1_\alpha(x_1,y_1)\otimes\mu^2_\beta(x_2,y_2),
\end{eqnarray*}
for all $x_1,y_1\in \mathcal {A}_1$, $x_2,y_2\in \mathcal {A}_2$.
Then, we have
\begin{eqnarray*}
(\nu^1_\alpha\otimes\nu^2_\beta)(x_1\otimes x_2,y_1\otimes
y_2)&=&\nu^1_\alpha(x_1,y_1)\otimes\nu^2_\beta(x_2,
y_2)\\&=&\alpha\circ\nu^1(x_1,
y_1)\otimes\beta\circ\nu^2(y_1, y_2)\\
&=&\nu^1(\alpha(x_1),\alpha(y_1))\otimes\nu^2(\beta(x_2),\beta(y_2))\\
&=&(\nu^1\otimes\nu^2)\big(\alpha(x_1)\otimes\beta(x_2),\alpha(y_1)\otimes\beta(y_2)\big)\\
&=&(\nu^1\otimes\nu^2)\big(\alpha\otimes\beta(x_1\otimes
x_2),\alpha\otimes\beta(y_1\otimes y_2)\big).
\end{eqnarray*}
Similarly, we have
\begin{eqnarray*}
(\mu^1_\alpha\otimes\mu^2_\beta)(x_1\otimes x_2,y_1\otimes y_2)
=(\mu^1\otimes\mu^2)\big(\alpha\otimes\beta(x_1\otimes
x_2),\alpha\otimes\beta(y_1\otimes y_2)\big).
\end{eqnarray*}
Now from the discussions above, we obtain
\begin{eqnarray*}
\mu^1_{\alpha}\otimes\mu^2_\beta&=&(\alpha\otimes\beta)\circ(\mu^1\otimes\mu^2)=(\mu^1\otimes\mu^2)_{\alpha\otimes\beta},\\
\nu^1_{\alpha}\otimes\nu^2_\beta&=&(\alpha\otimes\beta)\circ(\nu^1\otimes\nu^2)=(\nu^1\otimes\nu^2)_{\alpha\otimes\beta}.
\end{eqnarray*}
 Then, according to  Proposition \ref{prop3}, we have the following result (see also
[\ref{Yau5}, Corollary 3.6]):
\begin{theo} With notations above. Let $\big(\mathcal {A}_1,\mu^1,\nu^1\big)$ and
$\big(\mathcal {A}_2,\mu^2,\nu^2\big)$ be two Novikov-Poisson
algebras. If $\alpha$ and $\beta$ are algebra homomorphisms of
$\big(\mathcal {A}_1,\mu^1,\nu^1\big)$ and $\big(\mathcal
{A}_2,\mu^2,\nu^2\big)$, respectively. Then
 $\big(\mathcal {A}_1\otimes \mathcal
{A}_2,\mu^1_\alpha\otimes\mu^2_\beta,
\nu^1_\alpha\otimes\nu^2_\beta,\alpha\otimes\beta\big)$ forms a
Hom-Novikov-Poisson algebra.
\end{theo}

 We extend in the following the Yau's Theorem (see \cite{Yau}) to the Hom-Novikov-Poisson algebra case.
 The following theorem gives a way to construct Hom-Novikov-Poisson
algebras from a commutative Hom-associative algebra along with a
derivation.
\begin{theo} \label{theorem2} Let $(\mathcal {A},\cdot,\alpha)$ be a commutative
Hom-associative algebra and $\partial$ be a derivation such that
$\partial\alpha=\alpha\partial$. Define a new operation on
$\mathcal{A}$ by
\begin{eqnarray}\label{Gd-1}
x\star y=x\cdot \partial(y), \hspace{0.3cm} \mbox{for all} \
x,y\in\mathcal {A}.
\end{eqnarray}
Then $(\mathcal {A},\cdot,\star,\alpha)$ forms a Hom-Novikov-Poisson
algebra.
\end{theo}
\noindent{\it Proof.~}It follows from [\ref{Yau}, Theorem 1.2] that
$(\mathcal {A},\star,\alpha)$ is a Hom-Novikov algebra. We need to
check the two compatible conditions (\ref{Hom-NP-1}) and
(\ref{Hom-NP-2}). For any $x,y,z\in \mathcal {A}$, we have
\begin{eqnarray*}
(x\cdot y)\star \alpha(z)&=&(x\cdot y)\cdot\partial (\alpha(z))\\
&=&(x\cdot y)\cdot\alpha (\partial(z))\\
&=&\alpha(x)\cdot(y\cdot\partial(z))\\
&=&\alpha(x)\cdot(y\star z),
\end{eqnarray*}
which proves equation (\ref{Hom-NP-1}). Furthermore, we have
\begin{eqnarray*}
(x\star y)\cdot\alpha(z)-\alpha(x)\star(y\cdot
z)&=&(x\cdot\partial(y))\cdot\alpha(z)-\alpha(x)\cdot\partial(y\cdot z)\\
&=&(x\cdot\partial(y))\cdot\alpha(z)-\alpha(x)\cdot(\partial(y)\cdot z+y\cdot\partial(z))\\
&=&\alpha(x)\cdot(\partial(y)\cdot z)-\alpha(x)\cdot\big(\partial(y)\cdot z+y\cdot\partial(z)\big)\\
&=&-\alpha(x)\cdot\big(y\cdot\partial(z)\big)\\
&=&-(x\cdot y)\cdot\alpha(\partial(z)).
\end{eqnarray*}
Similarly, we have
\begin{eqnarray*}
(y\star x)\cdot\alpha(z)-\alpha(y)\star(x\cdot z)=-(y\cdot
x)\cdot\alpha(\partial(z)).
\end{eqnarray*}
Since $(\mathcal {A},\cdot,\alpha)$ is commutative, we get
\begin{eqnarray*}
(x\star y)\cdot\alpha(z)-\alpha(x)\star(y\cdot z)=(y\star
x)\cdot\alpha(z)-\alpha(y)\star(x\cdot z),
\end{eqnarray*}
which completes the proof.\QED\vskip7pt

The following result is due to Xu \cite{Xu2}.
\begin{prop}\label{p5}
Let $(\mathcal {A},\cdot,\ast)$ be a Novikov-Poisson algebra, such
that $(\mathcal {A},\cdot)$ is a commutative associative algebra and
$(\mathcal {A},\ast)$ is a Novikov algebra. Suppose that $(\mathcal
{A},\cdot)$ contains an identity element $1$, that is, $1\cdot
x=x\cdot 1=x$ for all $x\in\mathcal{A}$. Define a map
$\partial:\mathcal{A}\longrightarrow \mathcal{A}$ by
\begin{eqnarray}\label{der}
\partial (x)=1\ast x-(1\ast 1)\cdot x, \hspace{0.3cm} \mbox{for all} \
x\in\mathcal{A}.
\end{eqnarray}
Then $\partial$ is a derivation of $(\mathcal{A},\cdot)$.
\end{prop}

In the following, we give some examples of Hom-Novikov-Poisson
algebras by using Theorem \ref{theorem2} and Proposition \ref{p5}.

\begin{exam}\label{ex3}{\rm
 Assume that $(\mathcal{A},\cdot,\ast)$ is a
Novikov-Poisson algebra, where $(\mathcal{A},\cdot)$ is a
commutative associative algebra with unity $1$ and
$(\mathcal{A},\ast)$ is a Novikov algebra. Equation (\ref{NP-1})
implies that
$$x\cdot(1\ast 1)=x\ast 1,\hspace{0.3cm}\mbox{for all}\  x\in\mathcal{A}.$$
Since $(\mathcal{A},\cdot)$ is commutative, we have
\begin{eqnarray}\label{der1}
\partial (x)=1\ast x-x\ast 1, \hspace{0.3cm} \mbox{for all}\
x\in\mathcal{A},
\end{eqnarray}
in which $\partial$ is defined by (\ref{der}). Let $\alpha$ be an
algebra homomorphism of $(\mathcal{A},\cdot,\ast)$.  Define a new
multiplication on $\mathcal{A}$ by
\begin{eqnarray*}
x\bullet y=\alpha(x\cdot y), \hspace{0.3cm} \mbox{ for all} \
x,y\in\mathcal{A}.
\end{eqnarray*}
Then $(\mathcal{A},\bullet,\alpha)$ is a commutative Hom-associative
algebra. Obviously, $\alpha(1)=1$. As a consequence, $\alpha$
commutes with $\partial$. Furthermore, $\partial$ is a derivation of
$(\mathcal{A},\bullet)$ since it is a derivation of
$(\mathcal{A},\cdot)$ by Proposition \ref{p5}. Now define a new
operation on $\mathcal{A}$ by
$$x\star y=x\bullet\partial(y),\hspace{0.3cm}\mbox{for all}\ x,y\in\mathcal{A}.$$
Then, thanks to Theorem \ref{theorem2}, $(\mathcal {A},
\bullet,\star,\alpha)$ is a Hom-Novikov-Poisson algebra. }
\end{exam}
\begin{exam}\label{ex2}{\rm Let the ground field be the complex field
$\mathbb{C}$ and $\Delta$ be a nonzero abelian subgroup of
$\mathbb{C}$. Suppose that $f$ is a nontrivial homomorphism of
$\Delta$ into the additive group of $\mathbb{C}$. Consider a vector
space $\mathcal{A}$ with basis $\{x_a|a\in\Delta\}$. Fix an element
$q$ in $\Delta$. Define a multiplication on $\mathcal{A}$ by
\begin{eqnarray*}
x_a\cdot x_b=x_{a+b+q}, \hspace{0.3cm}\mbox{for all} \ a,b\in\Delta.
\end{eqnarray*}
Then $(\mathcal{A}, \cdot)$ is a commutative associative  algebra.
Define a linear map
$\partial_1:\mathcal{A}\longrightarrow\mathcal{A}$ by
\begin{eqnarray*}
\partial_1(x_a)=f(a+q)x_{a}, \hspace{0.3cm}\mbox{for all} \ a\in\Delta.
\end{eqnarray*}
It is easy to check that $\partial_1$ is a derivation of
$(\mathcal{A}, \cdot)$. Then we obtain a Novikov-Poisson algebra
$(\mathcal{A},\cdot,\ast)$ with the operation $\ast$ defined by
\begin{eqnarray*}
x_a\ast x_b=x_a\cdot\partial_1 (x_b)=f(b+q)x_{a+b+q},
\hspace{0.3cm}\mbox{for all} \ a,b\in\Delta.
\end{eqnarray*}
Note that $x_{-q}$ is an identity element of $(\mathcal{A},\cdot)$.
Define another linear map
$\partial_2:\mathcal{A}\longrightarrow\mathcal{A}$ by
\begin{eqnarray*}
\partial_2(x_a)= x_{-q}\ast x_a-x_{a}\ast x_{-q}, \hspace{0.3cm}\mbox{for all} \
a\in\Delta.
\end{eqnarray*}
Then it follows from Proposition \ref{p5} and the discussions
presented in Example \ref{ex3} that $\partial_2$ is a derivation of
$(\mathcal{A},\cdot)$.

Let $\alpha:\mathcal{A}\longrightarrow\mathcal{A}$ be a linear map
defined by
\begin{eqnarray*}
\alpha(x_a)=e^{a+q}x_{a}, \hspace{0.3cm}\mbox{for all} \ a\in\Delta.
\end{eqnarray*}
Then $\alpha$ is an algebra homomorphism of
$(\mathcal{A},\cdot,\ast)$ and $\alpha(x_{-q})=x_{-q}$. Moreover,
$\alpha$ commutes with $\partial_2$. Define a new multiplication on
$\mathcal{A}$ by
\begin{eqnarray*}
x_a\bullet x_b=\alpha(x_a\cdot x_b), \hspace{0.3cm}\mbox{for all} \
a,b\in\Delta,
\end{eqnarray*}
 which makes $(\mathcal{A},\bullet, \alpha)$ form a
commutative Hom-associative algebra. Also, $\partial_2$ is a
derivation of $(\mathcal{A},\bullet)$. Define a bilinear map $\star$
on $\mathcal{A}$ by
\begin{eqnarray*}
x_a\star x_b=x_a\bullet\partial_2( x_b), \hspace{0.3cm}\mbox{for
all} \ a,b\in\Delta,
\end{eqnarray*}
To be more precise, we have
\begin{eqnarray*}
x_a\star x_b=e^{a+b+2q}f(b+q)x_{a+b+q},\hspace{0.3cm}\mbox{for all}
\ a,b\in\Delta.
\end{eqnarray*}
According to Theorem \ref{theorem2}, $(\mathcal {A},
\bullet,\star,\alpha)$ is a Hom-Novikov-Poisson algebra.

}
\end{exam}
 \vskip12pt


\end{document}